\documentclass[11pt]{article}
\usepackage{amssymb, amsmath, amsfonts}

\begin{document}

\author{Cornelia Livia Bejan and Mircea Crasmareanu}
\title{Ricci solitons in manifolds with quasi-constant curvature}
\date{Dedicated to the memory of Stere Ianus (1939-2010)}
\maketitle

\begin{abstract}
The Eisenhart problem of finding parallel tensors treated already in the framework of quasi-constant curvature manifolds in \cite{x:j} is reconsidered for the symmetric case and the result is interpreted in terms of Ricci solitons. If the generator of the manifold provides a Ricci soliton then this is i) expanding on para-Sasakian spaces with constant scalar curvature and vanishing $D$-concircular tensor field and ii) shrinking on a class of orientable quasi-umbilical hypersurfaces of a real projective space=elliptic space form.
\end{abstract}



\medskip

\section*{Introduction}

In 1923, Eisenhart \cite{l:e} proved that if a positive definite
Riemannian ma\-ni\-fold $(M,g)$ admits a second order parallel
symmetric covariant tensor other than a constant multiple of the
metric tensor, then it is reducible. In 1926, Levy \cite{h:l}
proved that a parallel second order symmetric non-degenerated
tensor $\alpha $ in a space form is proportional to the metric
tensor. Note that this question can be considered as the dual to
the the problem of finding li\-near connections making parallel a
given tensor field; a problem which was considered by Wong in
\cite{y:w}. Also, the former question implies topological
restrictions namely if the (pseudo) Riemannian manifold $M$ admits
a pa\-ral\-lel symmetric $(0,2)$ tensor field then $M$ is locally the
direct product of a number of (pseudo) Riemannian manifolds,
\cite{h:w} (cited by \cite{g:z}). Another si\-tuation where the
parallelism of $\alpha $ is involved appears in the theory of
totally geodesic maps, namely, as is point out in \cite[p.
114]{o:c}, $\nabla \alpha =0$ is equivalent with the fact that
$1:(M,g)\rightarrow (M,\alpha )$ is a totally geodesic map.

While both Eisenhart and Levy work locally, Ramesh Sharma gives in \cite{r:s1} a global approach based on Ricci identities. In addition to space-forms, Sharma considered this {\it Eisenhart problem} in contact geometry \cite{r:s2}-\cite{r:s4}, for example for $K$-contact manifolds in \cite{r:s3}. Since then, several other studies appeared in various contact manifolds, see for example, the bibliography of \cite{c:c}.

Another framework was that of quasi-constant curvature in
\cite{x:j}; recall that the notion of {\it manifold with
quasi-constant curvature} was introduced by Bang-yen Chen and
Kentaro Yano in 1972, \cite{c:y}, and since then, was the subject
of several and very interesting works, \cite{b:r}, \cite{d:go},
\cite{yd:w}, in both local and global approaches. Unfortunately,
the paper of Jia contains some typos and we consider that a
careful study deserves a new paper. There are two remarks
regarding Jia result: i) it is in agreement with what happens in
all previously recalled contact geometries for the symmetric case,
ii) it is obtained in the same manner as in Sharma's paper
\cite{r:s1}. Our work improves the cited paper with a natural
condition imposed to the generator of the given manifold, namely
to be of torse-forming type with a regularity property.

Our main result is connected with the recent theory of Ricci
solitons \cite{ca:o}, a subject included in the Hamilton-Perelman approach
(and proof) of Poincar\'e Conjecture. A connection between Ricci
flow and quasi-constant curvature manifolds appears in \cite{c:z};
thus our treatment for Ricci solitons in quasi-constant curvature
manifolds seems to be new.

Our work is structured as follows. The first section is a very brief review of manifolds with quasi-constant curvature and Ricci solitons. The next section is devoted to the (symmetric case of) Eisenhart problem in our framework and the relationship with the Ricci solitons is pointed out. A technical conditions appears, which we call {\it regularity}, and is concerning with the non-vanishing of the Ricci curvature with respect to the generator of the given manifold. Let us remark that in the Jia's paper this condition is involved, but we present a characterization of these manifolds as well as some remarkable cases which are out of this condition namely: quasi-constant curvature locally symmetric and Ricci semi-symmetric metrics. A characterization of Ricci soliton is derived for dimension greater that 3.

Four concrete examples involved in possible Ricci solitons on quasi-constant manifolds are listed at the end. For the second example, we pointed out some consequences which are yielded by the hypothesis of compacity, used in paper \cite{d:g}, in connection with (classic by now) papers of T. Ivey and Perelman.

\section{Quasi-constant curvature manifolds. Ricci solitons}

Fix a triple $(M, g, \xi)$ with $M_n$ a smooth $n$($>2$)-dimensional manifold, $g$ a Riemannian metric on $M$ and $\xi $ an unitary vector field on $M$. Let $\eta $ the $1$-form dual to $\xi $ with respect to $g$.

If there exist two smooth functions $a$, $b\in C^{\infty }(M)$ such that:
$$
R(X, Y)Z=a[g(Y, Z)X-g(X, Z)Y]+b[g(Y, Z)\eta (X)-g(X, Z)\eta (Y)]\xi +
$$
$$
+b\eta (Z)[\eta (Y)X-\eta (X)Y] \eqno(1.1)
$$
then $(M, g, \xi )$ is called {\it manifold of quasi-constant
curvature} and $\xi $ is {\it the generator}, \cite{c:y}. Using the notation of \cite[p. 325]{d:t} let
us denote $M^n_{a, b}(\xi )$ this manifold.

It follows:
$$
R(X, Y)\xi = (a+b)[\eta (Y)X-\eta (X)Y] \eqno(1.2)
$$
$$
R(X, \xi )Z=(a+b)[\eta (Z)X-g(X, Z)\xi ] \eqno(1.3)
$$
while the Ricci curvature $S(X, Y)=Tr(Z\rightarrow R(Z, X)Y)$ is:
$$
S(X,Y)=[a(n-1)+b]g(X,Y)+b(n-2)\eta (X)\eta (Y) \eqno(1.4)
$$
which means that $(M, g, \xi)$ is an {\it eta-Einstein manifold}; in particular, if $a, b$ are scalars, then $(M, g, \xi )$ is an {\it quasi-Einstein manifold}, \cite{g:dt}. The scalar curvature is:
$$
r=(n-1)(na+2b), \eqno(1.5)
$$
and we derive:
$$
a=\frac{r-2S(\xi , \xi )}{(n-1)(n-2)}, \quad b=\frac{nS(\xi ,
\xi )-r}{(n-1)(n-2)}. \eqno(1.6)
$$
Then $a+b=\frac{S(\xi , \xi)}{n-1}$. Let us consider also the Ricci $(1, 1)$ tensor field $Q$
given by: $S(X, Y)=g(QX, Y)$. From $(1.4)$ we get:
$$
Q(X)=[a(n-1)+b]X+b(n-2)\eta (X)\xi \eqno(1.7)
$$
which yields:
$$
Q(\xi )=(a+b)(n-1)\xi \eqno(1.8)
$$
and then $\xi $ is an eigenvalue of $Q$.

In the last part of this section we recall the notion of Ricci solitons according to \cite[p. 139]{r:s5}. On the manifold $M$, a {\it Ricci soliton} is a triple $(g, V, \lambda )$ with $g$ a Riemannian metric, $V$ a vector field and $\lambda $ a real scalar such that:
$$
{\cal L}_Vg+2S+2\lambda g=0. \eqno(1.9)
$$
The Ricci soliton is said to be {\it shrinking, steady} or {\it expanding} according as $\lambda $ is negative, zero or positive.

Also, we adopt the notion of $\eta $-{\it Ricci soliton} introduced in the paper \cite{c:k} as a data $(g, V, \lambda , \mu)$:
$$
{\cal L}_Vg+2S+2\lambda g+2\mu \eta \otimes \eta =0. \eqno(1.10)
$$

\section{Parallel symmetric second order tensors and Ricci solitons}

Fix $\alpha $ a symmetric tensor field of $(0,2)$-type which we suppose to be parallel with respect to the Levi-Civita connection $\nabla $ i.e. $\nabla \alpha =0$. Applying the Ricci identity $\nabla ^2\alpha (X,Y;Z,W)-\nabla ^2\alpha (X,Y;W,Z)=0$ we obtain the relation $(1.1)$ of \cite[p. 787]{r:s1}:
$$
\alpha (R(X,Y)Z,W) +\alpha (Z,R(X,Y)W)=0 \eqno(2.1)
$$
which is fundamental in all papers treating this subject. Replacing $Z=W=\xi $ and using $(1.2)$ it results, by the symmetry of $\alpha $:
$$
(a+b)[\eta (Y)\alpha (X, \xi )-\eta (X)\alpha (Y, \xi )]=0. \eqno(2.2)
$$

\smallskip

{\bf Definition 2.1} $M^n_{a, b}(\xi )$ is called {\it regular} if $a+b\neq 0$.

\medskip

In order to obtain a characterization of such manifolds we
consider:

\medskip

{\bf Definition 2.2}(\cite{r:m}) $\xi $ is called {\it semi-torse forming
vector field} for $(M, g)$ if, for all vector fields $X$:
$$
R(X, \xi )\xi =0. \eqno(2.3)
$$

\smallskip

From $(1.2)$ we get: $R(X, \xi )\xi =(a+b)(X-\eta (X)\xi )$ and
therefore, if $X\in ker \eta=\xi ^{\bot}$, then $R(X, \xi )\xi
=(a+b)X$ and we obtain:

\medskip

{\bf Proposition 2.3} {\it For $M^n_{a, b}(\xi )$ the following are equivalent}: \\
i) {\it is regular}, \\
ii) $\xi $ {\it is not semi-torse forming}, \\
iii) $S(\xi , \xi )\neq 0$ {\it i.e.} $\xi $ {\it is non-degenerate with respect to} $S$, \\
iv) $Q(\xi )\neq 0$ {\it i.e. $\xi $ does not belong to the kernel of} $Q$. \\
{\it In particular, if $\xi $ is parallel} ($\nabla \xi =0$) {\it then
$M$ is not regular}.

\medskip

{\bf Remarks 2.4} i) From Theorems 2 and 3 of \cite[p. 175]{yd:w} a regular $M^n_{a, b}(\xi )$ is neither recurrent nor locally symmetric.\\
ii) From Theorem 3 of \cite[p. 228]{d:g} a regular $M^n_{a, b}(\xi )$ with $a$ and $b$ constants is not Ricci semi-symmetric.

\medskip

In the following we restrict to the regular case. Returning to $(2.2)$, with $X=\xi $ in:
$$
\eta (Y)\alpha (X, \xi )=\eta (X)\alpha (Y, \xi ) \eqno(2.4)
$$
we derive:
$$
\alpha (Y,\xi )=\eta (Y)\alpha (\xi , \xi)=\alpha (\xi , \xi)g(Y,
\xi). \eqno(2.5)
$$

The parallelism of $\alpha $ implies also that
$\alpha (\xi ,\xi )$ is a constant:
$$
X(\alpha (\xi , \xi ))=2\alpha (\nabla _X\xi , \xi)=2\alpha (\xi , \xi)g(\nabla _X\xi , \xi )=2\alpha (\xi , \xi )\cdot 0=0. \eqno(2.6)
$$

Making $Y=\xi $ in $(2.1)$ and using $(1.3)$ we get:
$$
\eta (Z)\alpha (X, W)-g(X, Z)\alpha (\xi , W)+\eta (W)\alpha (X, Z)-g(X, W)\alpha (\xi , Z)=0
$$
which yield, via $(2.5)$ and $W=\xi $:
$$
\alpha (X, Z) =\alpha (\xi , \xi )g(X, Z). \eqno(2.7)
$$
In conclusion:

\medskip

{\bf Theorem 2.5} {\it A parallel second order symmetric covariant
tensor in a regular $M^n_{a, b}(\xi )$ is a constant multiple of the metric tensor}.

\medskip

At the end of this section we include some applications of the above Theorem to Ricci solitons:

\medskip

Naturally, two remarkable situations appear regarding the vector field $V$: $V\in span{\xi }$ or $V\bot \xi$ but the second class seems far too complex to analyse in practice. For this reason it is appropriate to investigate only the case $V=\xi $. So,
we can apply the previous result for $\alpha:={\cal L}_{\xi }g+2S$ which yields $\lambda =-S(\xi , \xi )$.

\medskip

{\bf Theorem 2.6} {\it Fix a regular $M^n_{a, b}(\xi )$}.\\
i) {\it A Ricci soliton $(g, \xi , -S(\xi , \xi )\neq 0)$ can not be steady but is shrinking if the constant $S(\xi , \xi)$ is positive or expanding if $S(\xi , \xi)<0$}. \\
ii) {\it An $\eta $-Ricci soliton} $(g, \xi , \lambda , \mu)$ {\it provided by the parallelism of $\alpha +2\mu \eta \otimes \eta$ is given by}:
$$
\lambda +\mu =-S(\xi, \xi)\neq 0. \eqno(2.8)
$$
iii) {\it If $n\geq 4$ and $b\neq 0$ then $(g, \xi , -S(\xi , \xi ))$ is a Ricci soliton if and only if $\xi $ is geodesic i. e.} $\nabla _{\xi }\xi =0$ {\it and}:
$$
\frac{\xi (a+b)}{4b}+a(n-1)+b=\frac{a+b}{n-1}. \eqno(2.9)
$$

\medskip

{\bf Proof} iii) We have three cases:\\
I) $\alpha +2\lambda g=0$ on $span{\xi }$ yields the above expression of $\lambda $. \\
II) $\alpha +2\lambda g=0$ on $ker \eta=\xi ^{\bot }$ gives:
$$
\frac{\xi (a+b)}{4b}+\lambda +a(n-1)+b=0 \eqno(2.10)
$$
where we use the formula $(3.5)$ of \cite[p. 123]{g:m}. \\
III) $\alpha +2\lambda g=0$ on $(U, \xi)\in ker \eta \oplus span{\xi }$ gives:
$$
g(\nabla _U\xi , \xi)+g(U, \nabla _{\xi }\xi )=0.
$$
But the first term is zero since $\xi $ is unitary while the second implies that $\nabla _{\xi }\xi \in span{\xi }$. But again, $\xi $ being unitary we have that $\nabla _{\xi }\xi $ is orthogonal to $\xi $. \quad $\Box $

\medskip

{\bf Example 2.7} A para-Sasakian manifold with constant scalar curvature and vanishing $D$-concircular tensor is an $M^n_{a, b}(\xi )$ with \cite[p. 186]{d:g}:
$$
a=\frac{r+2(n-1)}{(n-1)(n-2)}, \quad b=\frac{-r-n(n-1)}{(n-1)(n-2)}
$$
and then, a Ricci soliton $(g, \xi )$ on it is expanding. This result can be considered as a version in para-contact geometry of the Corollary of \cite[p. 140]{r:s5} which states that a Ricci soliton $g$ of a compact $K$-contact manifold is Einstein, Sasakian and shrinking.

\medskip

From $(2.9)$ we get $r=-n$ and returning to formulae above it results:
$$
a=\frac{1}{n-1}, \quad b=\frac{-n}{n-1}.
$$

\medskip

{\bf Example 2.8} Let $N_{n+1}(c)$ be a space form with the metric $g$ and $M$ a {\it quasi-umbilical hypersurface} in $N$, \cite{c:y}, \cite[p. 175]{yd:w},
i.e. there exist two smooth functions $\alpha , \beta $ on $M$ and a 1-form $\eta $ of norm 1 such that the second fundamental form is:
$$
h_{ij}=\alpha g_{ij}+\beta \eta _i\eta _j.
$$
According to the cited papers $M$ is an $M^n_{a, b}(\xi )$ with:
$$
a=c+\alpha ^2, \quad b=\alpha \beta
$$
and $\xi $ the $g$-dual of $\eta $. This $M^n_{a, b}(\xi )$ is regular if and only if $c+\alpha ^2+\alpha \beta \neq 0$. Therefore, a Ricci soliton $(g, \xi )$ on this $M^n_{a, b}(\xi )$ is shrinking if $c+\alpha ^2+\alpha \beta >0$ and expanding if $c+\alpha ^2+\alpha \beta <0$.

\medskip

Inspired by Theorem 3 of \cite[p. 185]{d:g} let $N=\mathbb{R}P^{n+1}(c), c>0$ and $M$ an orientable quasi-umbilical hypersurface with $b=\alpha \beta >0$. Then:\\
i) a Ricci soliton $(g, \xi )$ on it is shrinking and $M$ is a real homology sphere (all Betti numbers vanish) if it is also compact, \\
ii) using the result of Ivey \cite{t:i}, for $n=3$ the manifold is of constant curvature being compact; so the case $n=4$ is the first important in any conditions or the case $n=3$ without compactness when we (possible) give up at the topology of real homology sphere, \\
iii) using again a classic result, now due to Perelman \cite{g:p}, the compactness implies that the Ricci soliton is gradient i.e. $\eta $ is exact.

\medskip

{\bf Example 2.9} Let $(M^{2n}_0, \omega _0, B)$ be a generalized Hopf manifold, \cite{d:o}, and $M^n$ an $n$-dimensional anti-invariant and totally geodesic submanifold. We set $\|\omega _0\|=2c$ and suppose that $B$ is unitary. Then, formula $(12.40)$ of \cite[p. 162]{d:o} gives that if $R^{\bot }=0$ then $M^n$ is of quasi-constant curvature with $a=c^2$ and $b=-\frac{1}{4}$. Therefore, $M^n$ is regular for $\|\omega _0\|\neq 1 $ and a Ricci soliton is shrinking if $\|\omega _0\|>1$ and expanding if $\|\omega _0\|<1$.

\medskip

{\bf Example 2.10} Suppose that $\xi $ is a {\it torse-forming vector field} i.e. there exist a smooth function $f$ and a $1$-form $\omega $ such that:
$$
\nabla _X\xi=fX+\omega (X)\xi. \eqno(2.11)
$$
From the fact that $\xi $ has unitary length it results $f+\omega (\xi )=0$ which means that $\xi $ is exactly a geodesic vector field.

\smallskip

{\bf Particular cases}: \\
i)(\cite{r:m}) If $\omega $ is exact then $\xi $ is called {\it concircular}; let $\omega =-du$ with $u$ a smooth function on $M$. Then $f=-\omega (\xi )=\xi (u)$. \\
ii) If $\omega =-f\eta $ then we call $\xi $ of {\it Kenmotsu type} since $(2.11)$ becomes similar to a expression well-known in Kenmotsu manifolds, \cite{c:c}.

\medskip

Let us restrict to ii). From $(2.11)$ a straightforward computation gives:
$$
R(X, Y)\xi =X(f)[Y-\eta (Y)\xi ]-Y(f)[X-\eta (X)\xi ] +f^2[\eta (X)Y-\eta (Y)X] \eqno(2.12)
$$
and a comparison with $(1.2)$ yields $a+b=-f^2$ and $f$ must be a constant, different from zero from regularity of the manifold. So, a possible Ricci soliton in a Kenmotsu type case must be expanding and with $S(\xi, \xi)$ and the scalar curvature constants, a result similar to Propositions 3 and 4 of \cite{c:c}.

\medskip

\noindent Seminarul Matematic "Al. Myller", \newline
 University "Al. I. Cuza"\newline
Ia\c si, 700506\newline
Romania\newline
e-mail: bejanliv@yahoo.com

\medskip

\noindent Faculty of Mathematics
\newline University "Al. I.Cuza" \newline
Ia\c si, 700506 \newline Romania \newline e-mail: mcrasm@uaic.ro

\smallskip

\noindent http://www.math.uaic.ro/$\sim$mcrasm

\end{document}